\documentclass{amsart}

\newtheorem{theorem}{Theorem}[section]
\newtheorem{lemma}[theorem]{Lemma}
\newtheorem{corollary}[theorem]{Corollary}

\theoremstyle{definition}

\theoremstyle{remark}

\numberwithin{equation}{section}

\usepackage{amssymb}
\usepackage{amscd}

\usepackage{color}  
\definecolor{darkgreen}{rgb}{0.03, 0.5, 0.03}
 \newcommand{\ec}{\color{black}} %

 \newcommand{\co} {\boldsymbol c}
 \newcommand{\Na} {\mbox{\rm Na}}
  \newcommand{\Max} {\mbox{\rm Max}}
    \newcommand{\Spec} {\mbox{\rm Spec}}
      \newcommand{\QMax} {\mbox{\rm QMax}}
        \newcommand{\QSpec} {\mbox{\rm QSpec}}

\newcommand{\fbar}{\boldsymbol{\overline{F}}}

  \newcommand{\stf} {\star{_{\!{_f}}}}
    \newcommand{\stt} {\widetilde{\star}}

\begin{document}

\title[Uppers to zero and semistar operations]
{Uppers to zero and semistar operations\\ in polynomial rings}

\author{Gyu Whan Chang}
\address{(Chang) Department of Mathematics, University of Incheon,
Incheon 402-749, Korea.} \email{whan@incheon.ac.kr}

\author{Marco Fontana}
\address{(Fontana) Dipartimento di Matematica, Universit\`a degli Studi  ``Roma Tre'', Largo San Leonardo Murialdo, 1 -- 00146 Roma, Italia.}
\email{fontana@mat.uniroma3.it}
\thanks{During the preparation of this paper, the second named author
was partially supported by MIUR, under Grant PRIN 2005-015278. }


\date{January 15 2007.}



\begin{abstract}  Given a stable semistar operation  of finite type  $\star$ on an integral domain $D$, we show that it   is   possible to define in a canonical way a stable semistar operation  of finite type  $[\star]$ on the polynomial ring $D[X]$, such that $D$ is a $\star$-quasi-Pr\"ufer domain if and only if each upper to zero in $D[X]$ is a quasi-$[\star]$-maximal ideal.  This result completes the investigation initiated by Houston-Malik-Mott \cite[Section 2]{hmm} in the star operation setting.   Moreover, we show that  $D$ is a Pr\"ufer $\star$-multiplication  (resp., a $\star$-Noetherian;  a $\star$-Dedekind) domain if and only if $D[X]$ is a Pr\"ufer $[\star]$-multiplication (resp.,   a $[\star]$-Noetherian; a $[\star]$-Dedekind) domain.  As an application of the techniques introduced here,  we obtain a  new interpretation of the  Gabriel-Popescu localizing systems of finite type on an integral domain $D$ (Problem 45 of \cite{cg}),  in terms of multiplicatively closed sets of the polynomial ring $D[X]$. 

\end{abstract}

\maketitle

\section{Introduction  and background results}


 Let $D$ be an integral domain with quotient field $K$. Let
$\boldsymbol{\overline{F}}(D)$ denote the set of all nonzero
$D$--submodules of $K$ and let $\boldsymbol{F}(D)$  (resp.,   $\boldsymbol{f}(D)$) be the set of
all nonzero fractional  (resp.,  finitely generated fractional) ideals of $D$.

Following Okabe-Matsuda \cite{o-m}, a \emph{semistar operation} on $D$ is a map $\star:
\boldsymbol{\overline{F}}(D) \to \boldsymbol{\overline{F}}(D), E
\mapsto E^\star$,  such that, for all $x \in K$, $x \neq 0$, and
for all $E,F \in \boldsymbol{\overline{F}}(D)$,
(a) $(xE)^\star=xE^\star$;
(b)  $E
\subseteq F$ implies $E^\star \subseteq F^\star$;
(c) $E \subseteq E^\star$ and $E^{\star \star} :=
\left(E^\star \right)^\star=E^\star$.
A \emph{(semi)star operation} is a semistar operation that,
restricted to $\boldsymbol{F}(D)$,  is a star operation (in the
sense of \cite[Section 32]{gilmer}). It is easy to see that a
semistar operation $\star$ on $D$ is a (semi)star operation if and
only if $D^\star = D$.

 If $\star$ is a semistar operation on $D$, then we can
consider a map\ $\star_{\!_f}: \boldsymbol{\overline{F}}(D) \to
\boldsymbol{\overline{F}}(D)$ defined by $E^{\star_{\!_f}}:=\bigcup \{F^\star\mid \ F \in
\boldsymbol{f}(D)$ and $F \subseteq E\}$, for each $E \in
\boldsymbol{\overline{F}}(D)$.
It is easy to see that $\star_{\!_f}$ is a semistar
operation on $D$, called \emph{the semistar operation of finite
type associated to $\star$}.  
 A semistar
operation $\star$ is called a \emph{semistar operation of finite
type} if $\star=\star_{\!_f}$.  It is easy to see that
$(\star_{\!_f}\!)_{\!_f}=\star_{\!_f}$ (that is, $\star_{\!_f}$ is
of finite type).

If $\star_1$ and $\star_2$ are two semistar operations on $D$, we
say that $\star_1 \leq \star_2$ if $E^{\star_1} \subseteq
E^{\star_2}$, for each $E \in \fbar(D)$. 
Obviously, for each semistar operation $\star$ defined on $D$, we
have $\star_{\!_f} \leq \star$. Let $d_D$ (or, simply, $d$)  be the \it identity (semi)star operation on $D$,   \rm clearly $d \leq \star$, for all semistar operation $\star$ on $D$.

We say that a nonzero ideal $I$ of $D$ is a
\emph{quasi-$\star$-ideal} if $I^\star \cap D = I$, a
\emph{quasi-$\star$-prime} if it is a prime quasi-$\star$-ideal,
and a \emph{quasi-$\star$-maximal} if it is maximal in the set of
all  proper  quasi-$\star$-ideals. A quasi-$\star$-maximal ideal is  a
prime ideal. It is possible  to prove that each  proper  quasi-$\star_{_{\!
f}}$-ideal is contained in a quasi-$\star_{_{\! f}}$-maximal
ideal.  More details can be found in \cite[page 4781]{fl}. We
will denote by $\QMax^{\star}(D)$  (resp., $\QSpec^\star(D)$) the set of the 
quasi-$\star$-maximal ideals  (resp., quasi-$\star$-prime ideals) of $D$.
When $\star$ is a (semi)star operation the notion of  quasi-$\star$-ideal coincides  with the ``classical'' notion of  \it  $\star$-ideal \rm (i.e., a nonzero ideal $I$ such that $I^\star = I$). 
%

  If $\Delta$ is a  nonempty  set of prime ideals of an integral
  domain
  $D$,  then the semistar operation $\star_\Delta$ 
  on $D$ defined by 
%
  $
  E^{\star_\Delta} := \bigcap \{ED_P \;|\;\, P \in \Delta\}, $
 for each $E \in \boldsymbol{\overline{F}}(D),
  $
   is called \it the spectral semistar operation associated to
  \rm
  $\Delta$.
  A semistar operation $\star$   on  an integral domain $D$ is
  called
  \it a
  spectral semistar operation \rm if there exists a nonempty subset $
  \Delta$ of the prime spectrum of $D$, $\mbox{\rm Spec}(D)$,  such that $\,\star =
  \star_\Delta$.

When $\Delta := \QMax^{\star_{_{\! f}}}(D)$, we set $\stt:= \star_{\Delta}$, i.e.,
%
  $E^{\stt} := \bigcap \left \{ED_P \mid P \in  \QMax^{\star_{_{\! f}}}(D) \right\}$,   for each $E \in \boldsymbol{\overline{F}}(D)$.
A semistar operation $\star$ is \emph{stable} if $(E \cap F)^\star
= E^\star \cap F^\star$, for each $E,F \in \fbar(D)$.
Spectral semistar operations are stable \cite[Lemma~4.1 (3)]{FH2000}. In particular, $\stt$ is a semistar operation stable and of finite type \cite[Corollary 3.9]{FH2000}.

 ~By $v_D$ (or, simply, by $v$) we denote  the $v$--(semi)star
operation defined as usual by  $E^v := (D:(D:E))$, for each $E\in
\boldsymbol{\overline{F}}(D)$. By  $t_D$ (or, simply, by $t$) we
denote  $(v_D)_{_{\! f}}$ the $t$--(semi)star operation on $D$ and
by  $w_D$ (or just by $w$) the stable semistar operation of finite
type associated to $v_D$ (or, equivalently, to $t_D$), considered
by F.G.  Wang  and R.L. McCasland in \cite{WMc97};  i.e. $w_D :=
\widetilde{v_D} = \widetilde{t_D}$.  Clearly $w_D\leq t_D \leq v_D$.  Moreover, it is easy to see that for each   (semi)star operation $\star$ of $D$, we have $\star \leq v_D$ and $\stf \leq t_D$ (cf. also \cite[Theorem 34.1 (4)]{gilmer}).


%


 Let $R$ be an overring of an integral domain $D$, let $\iota: D \hookrightarrow R$ be the canonical embedding 
and  let $\star$ be a semistar operation   on  $D$. We denote by  $\star_\iota$   the semistar operation  on   $R$ defined by $E^{\star_\iota} := E^\star$,  for each $E \in \fbar(R) \ (\subseteq \fbar(D))$.  
It is not difficult to see that 
  if $\star$ is a semistar operation of finite type (resp., a stable semistar operation)  on  $D$ then ${\star_\iota}$ is a semistar operation of finite type (resp., a stable semistar operation)  on   $R$ (cf. for instance \cite[Proposition 2.8]{fl2} and \cite[Propositions 2.11 and 2.13]{pi}).
  

A different approach to the stable semistar operation is possible by using the notion of localizing  system \cite{FH2000}.  Recall that a {\it
localizing system of ideals\/} ${\mathcal F}$ of $D$ is a set of (integral) ideals of $D$
verifying the following conditions
(a) if $
I \in\mathcal F  \mbox{ and if }  I\subseteq J $, then  $J\in\mathcal F
$;
(b) if 
$I\in\mathcal F$ and if  $J$ is an ideal of $D$ such that  $(J:_D iD)\in\mathcal F$,  for each $
i\in I$, then $ J\in\mathcal F$.
To avoid uninteresting cases, we assue that $\mathcal F$ is {\it nontrivial}, i.e., $\mathcal F$ is not empty and $(0) \not\in \mathcal F$.

The localizing systems, and the equivalent notions of Gabriel topologies (or, topologizing systems) and hereditary torsion theories, were introduced  in the 60's of the last century for the purpose of extending to non-commutative rings the theory of localization and for characterizing, from an ideal-theoretic point of view, the topologies associated to the hereditary torsion theories (cf. \cite{ga}, \cite[Ch. II, \S 2, Exercises 17-25, p. 157]{bk}, \cite{p},  and \cite[Ch. VI]{st}). 

For each  nonempty  subset $\Delta$ of prime ideals of $D$, set $\mathcal F(\Delta):=\{I \mbox{ ideal of }  D \mid I\not\subseteq P \mbox{ for
each } P\in\Delta\}\,.$
It is easy to verify that $\mathcal F(\Delta)$ is a localizing system of $D$ \cite[Proposition 5.1.4]{fhp}. If $P$ is a prime ideal of $D$, we denote
simply by $\mathcal F(P)$ the localizing system $\mathcal F (\{P\})$.
 It is obvious
that $\mathcal F(\Delta)=\bigcap\{\mathcal F(P)\mid P\in\Delta\}\ .$

\noindent \it A spectral localizing system \rm is a localizing system $\mathcal F$  such that  $\mathcal F = \mathcal F (\Delta)$, for some subset $\Delta $ of $\Spec(D)$.  A {\it localizing system of finite type\/} is a localizing system $\mathcal F$ 
such that for each $I\in\mathcal F$ there exists a finitely generated ideal
$J\in\mathcal F$ with $J\subseteq I$.

 Let $\mathcal F$ be a localizing system of ideals of $D$.   It is easy to see that, if $I,J \in  \mathcal F$,  then $IJ \in \mathcal F$, thus $\mathcal F$ is a multiplicative system of ideals and, inside the field of quotients $K$ of $D$,  it is possible to consider { \it  the generalized ring of fractions of $D$ with respect to $\mathcal F$}, i.e., $
D_{\mathcal F} := \bigcup \{(D:I) \mid I \in \mathcal F \} = \{ z \in K \mid (D:_D zD) \in \mathcal F \}\,.$
It is easy to see that, for each $E\in \fbar(D)$, 
$
E_{\mathcal F} := \bigcup \{(E:I) \mid I \in \mathcal F \} = \{ z \in K \mid (E:_D zD) \in \mathcal F \}
$
belongs to $ \fbar(D_{\mathcal F}) \ (\subseteq \fbar(D))$. We collect in the following lemma the main properties of the localizing systems that we will need in the present paper (cf. \cite[Proposition 2.4, Proposition 2.8, Theorem 2.10 (B) and Corollary 2.11]{FH2000} and 
 \cite[(5.1e), Lemma 5.1.5 (2), Propositions 5.1.4,   5.1.7 ((1)$\Leftrightarrow$(4)) and 5.18]{fhp}).

\begin{lemma} \label{F}   Let $\mathcal F$ be a localizing system   of ideals of   an integral domain $D$.
\begin{enumerate}
\item For each $E\in \fbar(D)$, the mapping $E \mapsto E_{\mathcal F}$ defines a stable semistar operation on $D$, denoted by $\star_{\mathcal F}$.
\item If $\boldsymbol{\Delta}(\mathcal F) := \{Q\in \Spec(D) \mid  Q\notin\mathcal F\}$, then $\mathcal F \subseteq \mathcal F(\boldsymbol{\Delta}(\mathcal F))$.

\item If $\mathcal F$ is a localizing system of finite type then $\mathcal F =  \mathcal F(\boldsymbol{\Delta}(\mathcal F))$.

\item If $\mathcal F = \mathcal F(\Delta)$ is a spectral localizing system then  $\mathcal F(\Delta) =  \mathcal F(\boldsymbol{\Delta}(\mathcal F))$.  Moreover, for each $E\in \fbar(D)$,
$
E_{\mathcal F(\Delta)} = \bigcap \{ED_P \mid P \in \Delta \}\,.$
\item   $\mathcal F$ is a localizing system of finite type if and only if there exists a quasi-compact subspace $\nabla$ of $\Spec(D)$ (endowed with the Zariski topology) such that $\mathcal F = \mathcal F(\nabla)$.

\item  Let  $\star$ be a semistar operation on $D$ and set $\mathcal F^\star := \{ I \mbox{ nonzero ideal of } D \mid I^\star = D^\star \}$. Then $\mathcal F^\star $ is a localizing system on $D$ and $\star_{\mathcal F^\star} = \star$ if and only if $\star$ is stable. 

\item The mapping $\mathcal F \mapsto \star_{\mathcal F}$ establishes a bijection between the set of the localizing systems (resp., the localizing systems of finite type) on $D$ and the set of the stable semistar operations (resp., the stable semistar operations of finite type) on $D$.

\end{enumerate}
\end{lemma}

The notion of quasi-Pr\"ufer domain (i.e.,  integral domain with Pr\"ufer integral closure) has a semistar operation analog introduced in \cite{cf}.
 The starting point of the present work is \cite[Corollary 2.4]{cf} where it is shown  that  the $t$-quasi-Pr\"ufer domains coincide with the UM$t$-domains (i.e., the integral domains such that each upper to zero in $D[X]$ is a maximal  $t_{D[X]}$-ideal).  There is no immediate extension to the semistar setting of the previous characterization, since  in the general case we do not have the possibility to work at the same time with a semistar  operation (like the $t$-operation) defined both on $D$ and on $D[X]$.  To overcome this difficulty,
given a semistar operation  of finite type  $\star$ on an integral domain $D$, we show that it  is  possible to define in a canonical way a semistar operation  of finite type  $[\star]$ on $D[X]$, such that $D$ is a $\star$-quasi-Pr\"ufer domain if and only if each upper to zero in $D[X]$ is a quasi-$[\star]$-maximal ideal. Moreover, we show that  $D$ is a P$\star$MD (resp., a $\stt$-Noetherian domain; a $\stt$-Dedekind domain) if and only if $D[X]$ is a P$[\star]$MD (resp.,   a $[\star]$-Noetherian domain; a $[\star]$-Dedekind domain).

As a by-product of the techniques introduced here, we obtain a  new interpretation of the  Gabriel  localizing   systems of finite type. More precisely,  we give an explicit natural bijection between the set of localizing systems of finite type $\mathcal F$ on an integral domain $D$ and the set of extended saturated multiplicative sets $\mathcal S$ of $ D[X]$; moreover,  $
E_{\mathcal F} = E\! \cdot\!D[X]_{{\mathcal{S}} } \cap K$, for all $E \in \fbar(D)$.

\section{Stable semistar operations and polynomial rings}

Let $D$ be an integral domain  with quotient field $K$, and let  $X$ be an
indeterminate over $K$.  For each polynomial $f \in K[X]$,  we denote by $\co_D(f)$ (or, simply, $\co(f)$) \it the content on $D$ of the polynomial $f$, \rm i.e., the  (fractional)  ideal of $D$ generated by the coefficients of $f$.
 

Let $\star$ be a semistar operation on $D$,  if     $\mathcal N^\star
:=  \{ g\in D[X] \mid  g \neq 0 \mbox{ and }  \co_D(g)^\star  =D^\star\}$,   then we set $\Na(D, \star) := D[X]_{\mathcal N^{\star}}$. The ring of rational functions $\Na(D, \star)$ is called \it  the $\star$--Nagata domain of $D$.  \rm When $\star =d$ the identity (semi)star operation on $D$,  then ${\mathcal N}^d ={\mathcal N}:=  \{ g\in D[X] \mid    \co_D(g) =D \}$. We set simply $\Na(D)$ instead of $\Na(D, d)= D[X]_{\mathcal N}$.  Note that  $\Na(D)$ coincides with the classical Nagata domain $D(X)$ (cf. for instance \cite[Chapter I, \S 6 page 18]{n} and \cite[Section 33]{gilmer}).

Recall  from \cite[Propositions 3.1 and 3.4]{fl} that:
\begin{enumerate}
\item[(a)]  $\mathcal N^\star= \mathcal N^{\stf}= \mathcal N^{\stt}=  D[X] \setminus \bigcup \{P[X] \mid P
  \in \QMax^{\stf}(D) \}$
is a saturated multiplicatively closed subset of  $D[X]$.

 \item[(b)]    $\Na(D, \star) = \Na(D, \stf) =\Na(D, \stt) =\bigcap\{ D_P(X) \mid P
  \in \QMax^{\stf}(D) \}$.

 \item[(c)]   
  $ \QMax^{\stf}(D)= \{M \cap D \mid M \in  \Max(\Na(D, \star)) \}\,.$

  \end{enumerate}
  
Furthermore, the stable semistar operation of finite type $\stt$  on $D$, canonically  associated to $\star$,  has the  following representation: 
$$E^{\tilde{\star}} = E\!\cdot\!\Na(D,\star)
  \cap K\,,  \mbox{      for each } E \in
  {\overline{\boldsymbol{F}}}(D)\,.
  $$

\smallskip

 More generally, let $R$ be an overring of $D$. We say that   $R$  is \it $t$-linked to $(D, \star)$ \rm if, for each nonzero finitely generated ideal $I$ of $D$, $I^{\star}=D^\star$ implies $(IR)^{t_R}= R$  \cite[Section 3]{eBF}.
 It is known that $R$ is a $t$-linked overring to $(D, \star)$ if and only if 
$R = R^{\stt}$ \cite[Lemma 2.9]{cf}.

Let  $\iota:  D \hookrightarrow R$ be the canonical embedding of $D$ in its overring $R$. If $R$ is a $t$-linked overring to $(D, \star)$ then  $(\stt)_\iota$ is a stable (semi)star operation of finite type on $R$ and 
$$ 
E^{(\stt)_\iota} = E\!\cdot\!\Na(D, \star) \cap K = E\!\cdot\!D[X]_{\mathcal{N}^\star} \cap K \,, \:\: \mbox{ for each  }  E \in \overline{\boldsymbol{F}}(R) 
$$
(cf. \cite[Lemma 2.9 ((i)$\Leftrightarrow$(v))]{cf} and the last part of Section 1).

At this point,  given an arbitrary  multiplicative subset $\mathcal S$   of $D[X]$, it is natural to ask whether the map $E \mapsto ED[X]_{\mathcal{S}}\cap K$,\  defined for all $E \in \overline{\boldsymbol{F}}(D)$, gives rise to a semistar operation $\star$ on $D$ (having the properties that $D^\star = R$, where $R :=D[X]_{\mathcal{S}}\cap K$, and that $R$ is $t$-linked to $(D, \star)$).   A complete answer to this question is given next. First we need a definition. Set:
$$\mathcal S^{\sharp} := D[X] \setminus  \bigcup \{ P[X] \mid  P \in \Spec(D) \mbox{ and }  P[X] \cap \mathcal S = \emptyset \}.
$$
  It is clear that $\mathcal S^{\sharp}$ is a saturated multiplicative set of $D[X]$ and that  $\mathcal S^{\sharp}$ contains the saturation of $\mathcal S$, i.e. $ \mathcal S^{\sharp} \supseteq  \overline{\mathcal S} =  D[X] \setminus  \bigcup \{ Q \mid  Q \in \Spec(D[X]) \mbox{ and } Q \cap \mathcal S = \emptyset \} $. We will call $\mathcal S^{\sharp}$ \it the extended saturation of ${\mathcal{S}}$ in $D[X]$ \rm and a multiplicative set  $\mathcal S$ of $D[X]$ is called \it extended saturated \rm  if $\mathcal S=\mathcal S^{\sharp}$.  Set
  $$
  \boldsymbol{\Delta}:= \boldsymbol{\Delta}(\mathcal{S}) := \{ P \in \Spec(D)\mid    P[X] \cap \mathcal{S} = \emptyset \}\,;
  $$
  obviously,  $\boldsymbol{\Delta}(\mathcal{S}) = \boldsymbol{\Delta}(\mathcal{S}^{\sharp})$. Let  $\boldsymbol{\nabla}:= \boldsymbol{\nabla}(\mathcal{S})$ be the set of the maximal elements of $\boldsymbol{\Delta}(\mathcal{S}) $.
\ec 

\begin{theorem} \label{circle} Let $\mathcal S$ be a multiplicative subset of   the polynomial ring   $D[X]$ and set
$E^{\circlearrowleft_{\mathcal{S}}}:= ED[X]_{\mathcal{S}}\cap K$,\  for all $E \in \overline{\boldsymbol{F}}(D)$. Clearly  $E^{\circlearrowleft_{\mathcal{S}}}\in \overline{\boldsymbol{F}}(D)$ and $ ED[X]_{\mathcal{S}} = E^{\circlearrowleft_{\mathcal{S}}}D[X]_{\mathcal{S}}$,\  for all $E \in \overline{\boldsymbol{F}}(D)$. 

\begin{enumerate}
\item[(a)]  The mapping ${\circlearrowleft_{\mathcal{S}}}: \overline{\boldsymbol{F}}(D)  \rightarrow
\overline{\boldsymbol{F}}(D)$, $E \mapsto E^{\circlearrowleft_{\mathcal{S}}}$ defines  a semistar operation on
$D$.

\item[(b)] ${\circlearrowleft_{\mathcal{S}}}$ is stable and of finite type, i.e.,  ${\circlearrowleft_{\mathcal{S}}}  =  \widetilde{\ {\circlearrowleft_{\mathcal{S}}} }$.

\item[(c)]  The extended saturation $ \mathcal{S}^{\sharp}$ of  $\mathcal S$  coincides with $\mathcal{N}^{{\circlearrowleft_{\mathcal{S}}}}:= \{ g \in D[X] \mid  g \neq~0 \mbox{ and } \co_D(g)^{\circlearrowleft_{\mathcal{S}}} = D^{{\circlearrowleft_{\mathcal{S}}}} \} $ and ${\circlearrowleft_{\mathcal{S}}} ={\circlearrowleft_{\mathcal{S}^\sharp}}$.

\item[(d)]  If $\mathcal{S}$ is extended saturated then  $\Na(D, {\circlearrowleft_{\mathcal{S}}}) = D[X]_{\mathcal{S}}$.

\item[(e)]  $\QMax^{\circlearrowleft_{\mathcal{S}}}(D) = \boldsymbol{\nabla}(\mathcal{S})$. In particular,  ${\circlearrowleft_{\mathcal{S}}}$ coincides with the spectral semistar operation $\star_{\boldsymbol{\nabla}(\mathcal{S})}$, i.e.,
$$
 E^{\circlearrowleft_{\mathcal{S}}} = \bigcap \{ED_P \mid P \in \boldsymbol{\nabla}(\mathcal{S})\}\,, \;\;\; \mbox{for all $E \in \overline{\boldsymbol{F}}(D)$}\,.
 $$

 \item[(f)]  ${\circlearrowleft_{\mathcal{S}}}$ is a (semi)star operation on $D$ if and only if  $\mathcal{S}
\subseteq \mathcal{N}^{v_D}$ or, equivalently, if and only if $D = \bigcap \{ D_P \mid P \in \boldsymbol{\nabla}(\mathcal{S})\}$. 

\item[(g)] The map $\mathcal{S} \mapsto {\circlearrowleft_{\mathcal{S}}}$ establishes a 1-1 correspondence between the extended saturated multiplicative subsets of $D[X]$ (resp., extended saturated multiplicative subsets of $D[X]$ contained in $\mathcal{N}^{v_D}$) and the set of the stable semistar (resp., (semi)star) operations of finite type on $D$.

\item[(h)]  Let $\mathcal{S}$ be an extended  saturated multiplicative set of $D[X]$. Then $\Na(D, v_D) =D[X]_{\mathcal{S}}$ if and only if ${\mathcal{S}} = \mathcal{N}^{v_D}$. 

\item[(i)]  Let $R:= D^{\circlearrowleft_{\mathcal{S}}} $ and let $\iota: D \rightarrow R$ be the canonical embedding.  The overring $R$ is $t$-linked to $(D, {\circlearrowleft_{\mathcal{S}}})$ and  ${\mathcal{S}} \subseteq {\mathcal{N}}^{v_R}$ (i.e., $({\circlearrowleft_{\mathcal{S}}})_\iota$ is a (semi)star operation on $R$). Moreover $({\circlearrowleft_{\mathcal{S}}})_\iota = w_R$ if and only if  the extended saturation $
{\mathcal{S}}^{\sharp_R}$ of the multiplicative set $\mathcal{S}$ in $R[X]$ coincides with $\mathcal{N}^{v_R}$.

\end{enumerate}
\end{theorem}

\begin{proof}  For the simplicity of notation,  set $\ast := {\circlearrowleft_{\mathcal{S}}}$.  Since $E \subseteq E^{\ast} $ and $E^{\ast} = ED[X]_{\mathcal{S}}  \cap K \subseteq  ED[X]_{\mathcal{S}} $, then  $E^{\ast} D[X]_{\mathcal{S}} =ED[X]_{\mathcal{S}} $.

(a)  The proof is straightforward.

(b) It suffices to show that $E^* \subseteq E^{\widetilde{\ast}}$ for each $E
\in \overline{\boldsymbol{F}}(D)$. If $0\neq x \in E^*$, then there exist $0\neq f \in
ED[X]$ and $0\neq g \in \mathcal{S}$ such that $x = \frac{f}{g} \in K$. So $xg = f$, and
thus $x\co_D(g)= \co_D(f) \subseteq E$. Note that $\co_D(g)^* = D^*$, since  $ gD[X]_{\mathcal{S}}  \subseteq \co_D(g)D[X]_{\mathcal{S}}$  and $ gD[X]_{\mathcal{S}}   = D[X]_{\mathcal{S}} $. Therefore  $g \in \mathcal{N}^{\ast} = \{ h \in D[X] \mid h \neq 0 \mbox{ and } \co_D(h)^* = D^* \}$   and so   $x = \frac{f}{g} \in ED[X]_{\mathcal{N}^{\ast}} \cap K = E\!\cdot\!\Na(D, \ast) \cap K  = E^{\widetilde{\ast}}$.  

(c) We have already observed (in the proof of (b)) that $\mathcal{S} \subseteq \mathcal{N}^{\ast} $. Since the multiplicative  set $\mathcal{N}^{\ast}$ coincides with  $D[X] \setminus \bigcup \{P[X]  \mid P \in \QMax^{\ast}(D)\}$ \cite [Proposition 3.1 (2)]{fl},  then $\mathcal{N}^{\ast} $ is extended saturated and so $\mathcal{S}^{\sharp} \subseteq \mathcal{N}^{\ast} $.
If  $0\neq g \in D[X]$ and $g\in {\mathcal{N}^{\ast}} \setminus {\mathcal{S}^\sharp}$ then $g \in Q[X]$, for some prime ideal $Q \in \Spec(D) \setminus \QMax^{\ast}(D)$ and $Q[X] \cap {\mathcal{S}} = \emptyset$. 
Note  that $Q^\ast \cap D \neq D$, i.e. $Q^\ast \neq D^\ast$, since $QD[X]_{\mathcal{S}} \neq D[X]_{\mathcal{S}}$. 
Since $\ast$ is a semistar operation of finite type,  we can find a quasi-$\ast$-maximal ideal $P$ in $D$ contaning  $Q^\ast \cap D $ and hence also containing $Q$. Therefore $g\in P[X]$, contradicting the assumption that $g \in \mathcal{N}^{\ast}$.   Finally, using (b), we have  ${\circlearrowleft_{\mathcal{S}}}  =  \widetilde{\ {\circlearrowleft_{\mathcal{S}}} } = \widetilde{\ast} =  {\circlearrowleft_{\mathcal{N}^\ast}}  =
{\circlearrowleft_{\mathcal{S}^\sharp}}$.

(d) is a straightforward consequence of (c).

(e)  By \cite [Proposition 3.1 (5)]{fl} and by (c) we have $\QMax^{\circlearrowleft_{\mathcal{S}}}(D) = \{ M \cap D \mid M\in \Max(D[X]_{\mathcal{N}^{\ast}})\} = \boldsymbol{\nabla}(\mathcal{S})$.  The remaining statement follows from (b).

(f) Suppose that $\ast'$  is a (semi)star operation on $D$, and let $g
\in \mathcal{S}$. If $g \not\in \mathcal{N}^{v_D}$, then $\co_D(g)^{-1} \neq D$, we can
choose $x \in \co_D(g)^{-1} \setminus D$, so   $x= \frac{xg}{g} \in
D[X]_{ \mathcal{S}} \cap K =D^{\ast'}$.  Since  $D = D^{\ast'}$  by assumption, we reach a contradiction. Thus $g
\in \mathcal{N}^{v_D}$.   Conversely, assume ${ \mathcal{S}} \subseteq \mathcal{N}^{v_D}$, then  $D^{\ast^\prime} = D[X]_{\mathcal{S}} \cap K \subseteq \Na(D, v) \cap K = D^w = D$.    The second equivalence follows from (e).

(g)  Let $\star$ be a stable semistar operation of finite type on $D$. Then $\star = \star_\Delta$, where $\Delta := \QMax^\star(D)$. Set $\boldsymbol{\mathcal{S}}(\Delta) := D[X] \setminus \bigcup \{P[X]  \mid P \in \Delta\}$. Clearly, $\boldsymbol{\mathcal{S}}(\Delta) $ is an extended saturated multiplicative set of $D[X]$ and  $\boldsymbol{\nabla}(\boldsymbol{\mathcal{S}}(\Delta) ) = \Delta$.  Therefore 
${\circlearrowleft_{\boldsymbol{\mathcal{S}}(\Delta)}} = \star_\Delta = \star$.   We easily conclude by using (b), (c) and (f).

(h) is a straightforward consequence of (g).

(i)  A part of this statement is a consequence of (f) and (h), after remarking that $({\circlearrowleft_{\mathcal{S}}})_\iota$ is a (semi)star operation on $R$ ``of type $\circlearrowleft$''  (defined by a multiplicative set of $R[X]$).   The fact that  $R$ is $t$-linked to $(D, {\circlearrowleft_{\mathcal{S}}})$ is a consequence of (b) and of \cite[Lemma 2.9 (i)$\Leftrightarrow$(v))]{cf}.
\end{proof}

The previous theorem leads to a  new interpretation of the   localizing systems of finite type on an integral domain $D$  in terms of multiplicatively closed sets of the polynomial ring $D[X]$.

\begin{corollary} The map $\mathcal F \mapsto \mathcal{S}:= \boldsymbol{\mathcal{S}}(\mathcal{F}) := D[X] \setminus \{Q[X] \mid Q \in  \Spec(D) \mbox{ and } Q \not\in \mathcal{F} \}$ establishes a natural bijection between the set of localizing systems of finite type $\mathcal F$ on an integral domain $D$ and the set of extended saturated multiplicative sets $\mathcal S$ of $ D[X]$.  Moreover,
$
E_{\mathcal F} = E\! \cdot\!D[X]_{\boldsymbol{\mathcal{S}}(\mathcal{F}) } \cap K \ (= E ^{\circlearrowleft_{{\boldsymbol{\mathcal{S}}(\mathcal{F}) } }}), \mbox{ for all } E \in \fbar(D)\,.
$

\end{corollary}

\begin{proof}   Let $\boldsymbol{\Delta}(\mathcal F) := \{ Q\in \Spec(D) \mid Q \not\in \mathcal{F} \}$  and so $\boldsymbol{\mathcal{S}}(\mathcal{F}) := D[X] \setminus \{Q[X] \mid Q \in \boldsymbol{\Delta}(\mathcal F) \}$.  Conversely, given an  extended saturated multiplicative set $\mathcal S$ of $ D[X]$, consider the set  $\boldsymbol{\Delta}(\mathcal{S}) := \{ P \in \Spec(D)\mid    P[X] \cap \mathcal{S} = \emptyset \}$ and define $\boldsymbol{\mathcal F}(\mathcal S) := \bigcap \{\mathcal F (P) \mid P \in  \boldsymbol{\Delta}(\mathcal{S}) \}$, where $\mathcal F (P) := \{I \mid I \mbox{ is an ideal of } D, I \not\subseteq P \}$.
The map defined by $\mathcal F \mapsto \boldsymbol{\mathcal{S}}(\mathcal{F})$  is a bijection, having as inverse the map defined by $  \mathcal{S} \mapsto \boldsymbol{\mathcal F}(\mathcal S)$.  As a matter of fact, given a  localizing systems of finite type $\mathcal F$ on $D$, then $\boldsymbol{\Delta}(\boldsymbol{\mathcal{S}}(\mathcal{F}) ) =  \boldsymbol{\Delta}(\mathcal{F})$ and thus $ \mathcal F = \boldsymbol{\mathcal F}(\boldsymbol{\mathcal{S}}(\mathcal{F}))$, since  for a localizing system of finite type we have $\mathcal F = \bigcap \{\mathcal F (P) \mid P \in  \boldsymbol{\Delta}(\mathcal{F}) \}$ \cite[Lemma 5.1.5 (2)]{fhp}.  Conversely, given an extended saturated multiplicative set $\mathcal S$ of $ D[X]$, then it is easy to see that 
$\boldsymbol{\Delta}(\mathcal S) \subseteq \boldsymbol{\Delta}(\boldsymbol{\mathcal F}(\mathcal S)) $. On the other hand,  if $Q \in  \boldsymbol{\Delta}(\boldsymbol{\mathcal F}(\mathcal S)) $, then $Q \not\in \boldsymbol{\mathcal F}(\mathcal S)$ and so $Q \not\in \mathcal F (P)$, i.e., $Q \subseteq P$, for some $P \in \boldsymbol{\Delta}(\mathcal S) $,  hence $Q[X] \cap \mathcal S = \emptyset$, i.e., $Q \in \boldsymbol{\Delta}(\mathcal S) $.  From the fact that $\boldsymbol{\Delta}(\mathcal S) = \boldsymbol{\Delta}(\boldsymbol{\mathcal F}(\mathcal S)) $ we have
$\boldsymbol{\mathcal{S}}(\boldsymbol{\mathcal F}(\mathcal S)) = 
 D[X] \setminus \{P[X] \mid P \in \boldsymbol{\Delta}(\mathcal S) \} = \mathcal S^{\sharp} = \mathcal S$. 
 
By Lemma \ref{F} (7), the last statement follows by observing that $\mathcal F$ coincides with $\mathcal F^{\circlearrowleft} := \{   I \mbox{ nonzero ideal of  } D  \mid I^{\circlearrowleft_{{\boldsymbol{\mathcal{S}}(\mathcal{F}) } }} \cap D = D \}$. 
 \end{proof}

\medskip


\ec

The notion of quasi-Pr\"ufer domain has a semistar analog introduced in \cite{cf}.
 Recall that an integral domain $D$ is a \it $\star$-quasi-Pr\"ufer  domain \rm if  for each prime ideal  $Q$  in $D[X]$ such that $Q \subseteq P[X]$, for some $P\in \QSpec^{\star}(D)$, then  $Q = (Q\cap D)[X]$.

As motivated in \cite{cf}, the previous notion has particular interest in case of semistar operations of finite type. Note that the $d$-quasi-Pr\"ufer  domains coincide with  the  quasi-Pr\"ufer  domains \cite[Theorem 1.1]{cf}.
For $\star =v$, we have observed in \cite[Corollary 2.4 (b)]{cf}  that  the $t$-quasi-Pr\"ufer domains coincide with the UM$t$-domains, i.e., the domains such that each upper to zero in $D[X]$ is a maximal  $t_{D[X]}$-ideal..  There is no immediate extension to the semistar setting of the previous characterization, since  in the general case we do not have the possibility to work at the same time with a semistar  operation (like the $t$-operation) defined both on $D$ and on $D[X]$. 

 This motivated the following question posed in \cite{cf}: \sl Given a semistar operation  of finite type  $\star$ on $D$, is it possible to define in a canonical way a semistar operation  of finite type  $\star_{\!_{D[X]}}$ on $D[X]$, such that $D$ is a $\star$-quasi-Pr\"ufer domain if and only if each upper to zero in $D[X]$ is a quasi-$\star_{\!_{D[X]}}$-maximal ideal~\!?  \rm 
 
In the next theorem and in the subsequent corollary we give a satistactory answer to the previous question, using the techniques introduced in Theorem \ref{circle}.

\begin{theorem}
\label{[*]} Let $D$ be an integral domain with quotient field $K$, let $X, Y$ be two 
indeterminates over $D$ and let $\star$ be a semistar operation on $D$.   Set $D_1 := D[X]$, $K_1 := K(X)$ and take the following subset of $\Spec(D_1)$:
$$
\boldsymbol{\Delta}_1^\star :=  \{Q_1 \in \Spec(D_1) \mid Q_1 \cap D = (0) \mbox{ or } Q _1=
(Q_1 \cap D)[X] \mbox{ and } (Q_1 \cap D)^{\stf} \subsetneq D^{\star}\}\,.
$$
Set  $\mathcal{S}_1^{\star}:= \mathcal{S}(\boldsymbol{\Delta}_1^\star ) :=  D_1[Y] \setminus
\left(\bigcup\{Q_1[Y] \mid Q_1 \in \boldsymbol{\Delta}_1^\star \}\right)$ and:
$$
E^{\circlearrowleft_{\mathcal{S}_1^{\star}} }:= E[Y]_{\mathcal{S}_1^{\star}} \cap K_1\,, \;\;\;  \mbox{ for all } E \in
\fbar(D_1).
$$

\begin{enumerate}
\item[(a)]  The mapping  $[\star]:= {\circlearrowleft_{\mathcal{S}_1^{\star}} } :  \fbar(D[X]) \rightarrow \fbar(D[X])$, 
$E\mapsto E^{\circlearrowleft_{\mathcal{S}_1^{\star}} }$
is a stable semistar operation of finite type on $D[X]$, i.e., $\widetilde{\ [\star]\ } = [\star]$.  Moreover, if $\star$ is a (semi)star operation on $D$, then $[\star]$ is a (semi)star operation on $D[X]$.

\item[(b)] $[\ \! \stt\ \! ] = [\stf] =[\star]$.

\item[(c)]  $(ED[X])^{[\star]} \cap K = ED_1[Y]_{{\mathcal{S}_1^{\star}}} \cap K = E^{\stt}$ for all $E \in \fbar(D)$.

\item[(d)] $(ED[X])^{[\star]} = E^{\stt}D[X]$,  for all $E \in \fbar(D)$.

\item[(e)]  $\QMax^{[\star]}(D_1) = \{Q_1 \mid  Q_1 \in \Spec(D_1)$ such that $Q_1 \cap D = (0)$ and
$\co_D(Q_1)^{\stf} = D^\star \} \cup \{P[X] \mid  P \in  \QMax^{\stf}(D)\}$. 
\item[(f)]   $[w_D] = [t_D]  = [v_D]  =  \widetilde{v_{D_1}}= w_{D_1}$.
\end{enumerate}
\end{theorem}

\begin{proof}   Note that,  if  $Q_1 \in \Spec(D[X])$ is not an upper to zero  and $(Q_1 \cap D)^{\stf} \subsetneq D^{\star}$,  then  the prime ideal   $Q_1 \cap D$   is contained in a quasi-$\stf$-maximal ideal of $D$.   Moreover if $Q_1 \cap D = (0)$ and 	$\co_D(Q_1)^{\stf}  \subsetneq D^{\star}$ then $\co_D(Q_1)^{\stf}$ is contained in a quasi-$\stf$-prime ideal $P$ of $D$   and hence  $Q_1 \subseteq P[X]$ with $P^{\stf} \subsetneq D^{\star}$.  Set  $\boldsymbol{\nabla}_1^\star :=  \{Q_1 \in \Spec(D_1) \mid $ either  $Q_1 \cap D = (0)$ and  $\co_D(Q_1)^{\stf}  =D^\star$ or
$Q_1= PD[X] \mbox{ and } P \in \QMax^{\stf}(D)\}$.  

It is easy to see that $$
\mathcal{S}_1^{\star}:= D_1[Y] \setminus
\left(\bigcup\{Q_1[Y] \mid Q_1 \in \boldsymbol{\Delta}_1^\star \}\right) = D_1[Y] \setminus
\left(\bigcup\{Q_1[Y] \mid Q_1 \in \boldsymbol{\nabla}_1^\star \}\right) =  \mathcal{S}(\boldsymbol{\nabla}_1^\star ).
$$

(a) follows from Theorem  \ref{circle}  ((a), (b) and (f)).

(b) Since $\QMax^{\stf}(D) =  \QMax^{\stt}(D)$, the conclusion follows easily from  the fact that $ \mathcal{S}_1^{\stt} = \mathcal{S}_1^{\stf} = \mathcal{S}_1^\star$. 

(c) Let $\mathcal{N}^{[\star]} := \{ g \in D_1[Y] \mid g \neq 0 \mbox{ and }  \co_{D_1}(g)^{[\star]}= D_1^{[\star]}\}$. Since by construction  $\mathcal{S}_1^{\star}$ is an extended saturated multiplicative set of $D_1$ we know  that  $\mathcal{S}_1^{\star} = \mathcal{N}^{[\star]}$ (Theorem \ref{circle} (c)).  On the other hand, if $h\in \mathcal{N}^{\star} =  D[X]  \setminus \left(\bigcup \{ P[X] \mid P \in \QMax^{\stf}(D) \}\right)$  then  $h \in  D[X][Y] \setminus \left(\bigcup\{Q_1[Y] \mid Q_1 \in \boldsymbol{\nabla}_1^\star \} \right) = \mathcal{N}^{[\star]}$. Therefore, for all $E \in \fbar(D)$,  $E^{\stt} = ED[X]_{\mathcal{N}^{\star}} \cap K  \subseteq ED_1[Y]_{ \mathcal{N}^{[\star]}} \cap K = 
\left(ED_1[Y]_{ \mathcal{N}^{[\star]}} \cap K_1 \right) \cap K = 
(ED_1)^{[\star]} \cap K =  (ED[X])^{[\star]} \cap K $.

%


For the reverse containment, let
 $0\neq z = \frac{f}{g} \in ED[X][Y]_{\mathcal{S}_1^{\star}}  \cap K$,  where $ z \in K$ and
  $ f, g \in K[X, Y]$ are nonzero polynomials
  such that  $  f \in ED[X][Y]$  and $g \in
{\mathcal{S}_1^{\star}} = \mathcal{N}^{[\star]} $.  
  Set $g = g_0 + g_1Y + \cdots +g_nY^n$, where $g_i \in D_1$ with $g_n \neq 0$ and $n \geq 0$; then $c_{D_1}(g) = (g_0, g_1, \dots , g_n)$ and $c_D(g) = c_D(g_0) + c_D(g_1) + \cdots +c_D(g_n)$.   
Let $Q_1 \in \boldsymbol{\Delta}_1^\star$. Since $\co_{D_1}(g)^{[\star]}= D_1^{[\star]}$, then $g \not\in Q_1[Y]$, and hence $(g_0, g_1, \dots , g_n) \nsubseteq Q_1$. So at least one among the $g_i$'s  is not contained in $Q_1$, and thus  $\co_D(g) \nsubseteq  Q_1 \cap D$.  In particular  $\co_D(g) \nsubseteq P$, for all $P \in \QMax^{\stf}(D)$, i.e., $\co_D(g)^{\stt} = D^{\stt}$.  \ec On the other hand, 
$z\co_{D}(g)= \co_{D}(zg)= \co_{D}(f)  \subseteq E$. Hence  $z \in z\co_D(g)^{\stt}  \subseteq E^{\stt}$. 
 Therefore we conclude that $ED[X]_{{\mathcal{N}}^{\star} } \cap K =  E^{\stt}$.

(d) By (c),  $(ED[X])^{[\star]} \cap K = E^{\stt}$, and thus
$E^{\stt}D[X] \subseteq (ED[X])^{[\star]}$, for all $E \in \fbar(D)$.


For the converse, let
$0\neq \frac{h}{\ell} = \frac{f}{g} \in (ED[X])^{[\star]} = ED[X][Y]_{\mathcal{S}_1^\star} \cap
K_1$, where $h, \ell \in K[X]$ are nonzero polynomials such that GCD$(h, \ell ) = 1$ in
$K[X]$, $0\neq f \in ED[X][Y]$, and $0\neq g \in {\mathcal{S}_1^\star} $. Then $\ell f =
h g$, and since $K[X, Y]$ is a UFD and GCD$(h, \ell)=1$, we have
$\ell \mid g$ in $K[X, Y]$, i.e., $g= \ell\!\cdot\! \gamma $ for some $\gamma \in
 K[X,Y]$. 

We claim that $\ell  \in K$.  Assume that  $\ell  \in K[X]
\setminus K$.
Choose a prime ideal $Q_1$ of
$D_1= D[X]$ such that $\ell K[X] \cap D[X] \subseteq Q_1$ and $Q_1 \cap D =
(0)$. Then $g=  \ell\!\cdot\! \gamma \in Q_1K[X,Y] \cap D[X,Y] =  Q_1[Y] \in \boldsymbol{\Delta}^\star_1$ and so $g \not\in { \mathcal{S}}_1^\star$, which is a contradiction.

Since $0 \neq \ell \in K$, set  $h' := \frac{h}{\ell} \in K[X]$.
Then $h' = \frac{f}{g} \in ED[X][Y]_{\mathcal{S}_1^\star} \cap K[X]$
and so $h'\! \cdot\! g = f$.
Since $g \in \mathcal{S}_1^\star$, by the proof of (c) above, we have $\co_D(g)^{\stt} = D^{\stt}$,
and hence $\co_{D}(h') \subseteq \co_{D}(h')\co_D(g)^{\stt} \subseteq (\co_D(h')\co_D(g)^{\stt})^{\stt}
 = (\co_D(h')\co_D(g))^{\stt} = \co_{D}(h'\! \cdot\! g)^{\stt} =
\co_{D}(f)^{\stt}  \subseteq E^{\stt}$  (cf. \cite[Corollary 28.3]{gilmer} for the fourth equality).    We conclude that $h' =\frac{f}{g} \in E^{\stt}D[X]$.

(e)  By \cite[Proposition 3.1 (5)]{fl} we know that $\QMax^{[\star]}(D_1)  = \{ M \cap D_1 \mid M \in \Max(D_1[Y]_{\mathcal{N}^{[\star]}}) \} $  and it is easy to verify  that this last set coincides with $\boldsymbol{\nabla}_1^\star $.

%

(f) If 
$\stf= t$, then  by (e) $\QMax^{[\star]}(D[X]) = \{Q_1 \mid  Q_1\in \Spec(D_1)$ such that $Q_1 \cap D = (0)$ and
$\co_D(Q_1)^{t} = D\} \cup  \{P[X] \mid  P \in  \Max^{t}(D)\} 
 =  \Max^{t_{D[X]}}(D[X])$. The last equality holds because it is wellknown that if $P \in  \Max^{t}(D)$ then $P[X] \in \Max^{t_{D[X]}}(D[X])$ \cite[Proposition 4.3]{Hedstrom/Houston: 1980}  and \cite[Proposition 1.1]{hz2};  moreover, if $ Q_1\in \Spec(D_1)$ is such that $Q_1 \cap D = (0)$, then  $Q_1$ is a $t_{D_1}$-maximal ideal  if and only if  $\co_D(Q_1)^{t} =D$  \cite[Theorem 1.4]{hz2}. Thus, by (a) and (b) and by the fact that $\QMax^{[\star]}(D[X])= \Max^{t_{D[X]}}(D[X])$, we have  $[v_D] =[t_D] =[w_D] =\widetilde{\ [w_D]\ } = \widetilde{v_{D[X]}}= w_{D[X]}$.
\end{proof}

\begin{corollary}
Let $\star$ be a semistar on an integral domain $D$ and let $[\star]$ be the stable semistar operation of finite type on $D[X]$ canonically associated to $\star$ as in Theorem \ref{[*]} (a). The following statements are equivalent:
\begin{enumerate}
\item[(i)] $D$ is a
$\stf$-quasi-Pr\"ufer domain.
\item[(ii)] $D[X]$ is a
$[\star]$-quasi-Pr\"ufer domain.
\item[(iii)] Each upper to zero  is a quasi-$[\star]$-maximal ideal of $D[X]$.
\end{enumerate}
\end{corollary}

\begin{proof} The equivalence (i)$\Leftrightarrow$(iii) follows easily from Theorem \ref{[*]} (e) and from the fact that $D$ is a $\stf$-quasi-Pr\"ufer  domain if and only if,    for each  upper
to zero $Q$   in $D[X]$, $\co(Q)^{\stf}= D^{\star}$  \cite[Lemma 2.3]{cf}..

For the equivalence between (i) and (ii), recall that $D$ is a
$\stf$-quasi-Pr\"ufer domain if and only if  $D_P$ is a quasi-Pr\"ufer domain, for each  quasi-$\stf$-maximal ideal   $P$  of $D$  \cite[Theorem 2.16 ($(1_{\stf})\Leftrightarrow(11_{\stf})$)]{cf}.   Moreover,  for each prime ideal $P$ of $D$,  $D[X]_{P[X]}$ coincides with the Nagata ring $D_P(X)$  and this is a quasi-Pr\"ufer domain if and only if $D_P$ is a quasi-Pr\"ufer domain  \cite[Theorem 1.1((1) $\Leftrightarrow$ (9))]{cf}.   

(i)$\Rightarrow$(ii)  Since we know already that  (i)$\Leftrightarrow$(iii), in the present situation  
we have $\QMax^{[\star]}(D[X]) = \{Q_1 \mid  Q_1 \in \Spec(D_1)$ such that $Q_1 \cap D = (0) \} \cup \{P[X] \mid  P \in  \QMax^{\stf}(D)\}$.  The conclusion follows from the fact that  $D[X]_{Q_1}$ is clearly a DVR for each upper to zero $Q_1$ of $D[X]$ and $D[X]_{P[X]} = D_P(X)$ is a quasi-Pr\"ufer domain, since $D_P$  is a quasi-Pr\"ufer domain, for each  for each  quasi-$\stf$-maximal ideal   $P$  of $D$.

(ii)$\Rightarrow$(i) is obvious by the previous argument.
\end{proof}

From the previous corollary  and from \cite[Corollary 2.4 (b)]{cf}, we re-obtain   that an integral domain $D$ is a UM$t$-domain if and only if the polynomial ring $D[X]$ is a UM$t$-domain \cite [Theorem 2.4]{fgh}, since by  Theorem \ref{[*]} (f), the semistar operation $[t_D]$ on $D[X]$ coincides with $w_{D[X]}$  and the notions of $w$-quasi-Pr\"ufer domain and  $t$-quasi-Pr\"ufer domain coincide.

\smallskip

 Let $\star$ be a semistar operation on an integral domain $D$. We say that $D$ is a \it $\star$-Noetherian domain \rm  if $D$ has the ascending chain condition on quasi-$\star$-ideals of $D$. It is easy to show that $D$ is $\star$-Noetherian if and only if each nonzero ideal $I$ of $D$ is $\stf$-type, i.e., $I^{\stf} = J^{\stf}$ for some $J \in \boldsymbol{f}(D)$ and $J \subseteq I$.  It is known that $D$ is $\stt$-Noetherian if and only if $\Na(D, \star) = D[X]_{\mathcal N^{\star}}$ is Noetherian,  \cite[Theorem 4..36]{pi} (cf. \cite[Theorem 2.6]{chang} for the star operation case). An $I \in \boldsymbol{\overline{F}}(D)$ is said to be \it quasi-$\star$-invertible \rm (resp., {\it $\star$-invertible}) if $(I:(D^{\star}:I))^{\star} = D^{\star}$ (resp., $(I:(D:I))^{\star} = D^{\star}$).  Recall that $D$ is   a \it $\star$-Dedekind domain \rm   if each nonzero (integral) ideal of $D$ is quasi-$\stf$-invertible and $D$ is  a \it  Pr\"{u}fer $\star
$-multiplication domain  \rm (for short,  \it P$\star $MD\rm) if
every nonzero finitely generated (integral) ideal of  $D$ is $\star
_{_{\!f}}$-invertible (cf. for instance \cite{fjs}). It is known that $D$ is a $\star$-Dedekind domain  if and only if $D$ is a P$\star$MD and a $\star$-Noetherian domain \cite[Proposition 4.1]{eBFP}.   

\begin{corollary}
Let $\star$ be a semistar on an integral domain $D$ and let $[\star]$ be the stable semistar operation of finite type on $D[X]$ canonically associated to $\star$ as in Theorem \ref{[*]} (a). Then
\begin{enumerate}
\item $D$ is a P$\star$MD if and only if $D[X]$ is a P$[\star]$MD.
\item $D$ is a $\stt$-Noetherian domain if and only if $D[X]$ is a $[\star]$-Noetherian domain.
\item $D$ is a $\stt$-Dedekind domain  if and only if $D[X]$ is a $[\star]$-Dedekind domain.
\end{enumerate}
\end{corollary}

\begin{proof}
(1) By Theorem 2.3(d), we have $(D[X])^{[\star]} = D^{\stt}[X]$, and hence $(D[X])^{[\star]}$ is integrally closed if and only if $D^{\stt}$ is integrally closed. Thus the result follows directly from Corollary 2.4 and \cite[Corollary 2.17]{cf}.

(2) Assume that $D$ is a $\stt$-Noetherian domain.
Then $D[X]_{\mathcal N^{\star}}$ is Noetherian and so  $(D[X]_{\mathcal N^{\star}})[Y]= (D[X][Y])_{\mathcal N^{\star}}$ is also Noetherian.
On the other hand, recall that ${\mathcal{N}^{\star}} \subseteq \mathcal{N}^{[\star]}$ (cf. the proof of Theorem 2.3(c)), and so $(D[X][Y])_{\mathcal{N}^{[\star]}} = ((D[X][Y])_{\mathcal{N}^{\star}})_{\mathcal{N}^{[\star]}}$ is Noetherian. Hence,
 $D[X]$ is $[\star]$-Noetherian.

For the converse, let $I$ be a nonzero ideal of $D$. Since $D[X]$ is $[\star]$-Noetherian, then $(ID[X])^{[\star]} = (f_1, f_2, \dots , f_n)^{[\star]}$, for a finite family of polynomials $f_1, f_2 \dots , f_n \in ID[X]$.  Set $J = \co_D(f_1) +  \co_D(f_1)+ \cdots + \co_D(f_n)$.  Clearly $(f_1, f_2, \dots , f_n) \subseteq JD[X]$ and thus $(ID[X])^{[\star]} = (JD[X])^{[\star]}$. Therefore, by Theorem 2.3(c), we have $I^{\stt} = (ID[X])^{[\star]} \cap K = (JD[X])^{[\star]} \cap K = J^{\stt}$ and so we conclude that $D$ is $\stt$-Noetherian.

(3) This is an immediate consequence of (1), (2) and \cite[Proposition 4.1]{eBFP}.  
\end{proof}

\ec

\bibliographystyle{amsplain}

\begin{thebibliography}{10}

 





%

%

\bibitem{bk} N. Bourbaki, Alg\`ebre Commutative, Chapitres I-II, Hermann, Paris, 1961.

%

%

\bibitem{chang} G.W. Chang, {\em $*$-Noetherian domains and the
ring $D[X]_{N_*}$}, J. Algebra 297(2006), 216-233.



\bibitem{cf} G.W. Chang and M. Fontana, \it  Upper to zero in polynomial rings and Pr\"ufer-like domains.  \rm Preprint.

\bibitem{cg} S.T. Chapman and S. Glaz, {\em One hundred problems in commutative ring theory}, ``Non-Noetherian Commutative Ring Theory'' (S.~T. Chapman and
  S. Glaz, eds.), Kluwer Academic Publishers, 2000, pp.~459-476. 



   
%








 
\bibitem{eBF} S. El Baghdadi and M. Fontana, \it Semistar linkedness and flatness, Pr\"ufer semistar multiplication domains, \rm Comm. Algebra \bf 32 \rm (2004), 1101--1126.




\bibitem{eBFP} S. El Baghdadi,  M. Fontana, and G. Picozza,  \it Semistar Dedekind domains, \rm  J. Pure Appl. Algebra \bf 193 \rm  (2004), 27--60. 

  

   
\bibitem{FH2000}
M. Fontana and J.A. Huckaba, \emph{Localizing systems and
semistar
  operations}, ``Non-Noetherian Commutative Ring Theory'' (S.~T. Chapman and
  S. Glaz, eds.), Kluwer Academic Publishers, 2000, pp.~169-198.
   

\bibitem{fhp} M. Fontana, J. Huckaba, and I. Papick,  {Pr\"ufer
domains}, Marcel Dekker, 1997.


\bibitem{fjs} M. Fontana, P. Jara, and E. Santos, \it
Pr\"ufer $\star$-multiplication domains and semistar operations, \rm J. Algebra Appl. \bf 2 \rm (2003), 21-50.


\bibitem{fgh} M. Fontana, S. Gabelli, and E. Houston, {\em
UMT-domains and domains with Pr\"ufer integral closure}, Comm.
Algebra 26(1998), 1017-1039.

   

   

\bibitem{fl2} M. Fontana and K.A. Loper, {\em Kronecker function
rings: a general approach}, in Ideal Theoretic Methods in
Commutative Algebra, Lecture Notes in Pure Appl. Math., Marcel
Dekker, 220(2001), 189-205.

\bibitem{fl} M. Fontana and K.A. Loper, {\em Nagata rings,
Kronecker function rings and related semistar operations}, Comm.
Algebra 31(2003), 4775-4801.




   
   \bibitem{ga} P. Gabriel, {\em Des cat\'egorie ab\'eliennes}, Bull. Soc. Math. France 90(1962), 323-448.
   

\bibitem{gilmer} R. Gilmer, {Multiplicative Ideal Theory},
 Marcel Dekker,  New York, 1972.

 
\bibitem{GH}
R. Gilmer and J. F. Hoffmann, \it A characterization of Pr\"ufer domains in terms of polynomials, \rm Pacific J. Math., 60(1975), 81-85. 





%

   
%
   

%


   
%
%
%

\bibitem{Hedstrom/Houston: 1980} J.R. Hedstrom and E.G. Houston, {\em Some
remarks on star-operations}, J. Pure Appl. Algebra 18(1980), 37--44.


%
%

%




\bibitem{hous}  E. Houston, \it  Uppers to zero in polynomial rings, \rm
in ``Multiplicative Ideal Theory in Commutative Algebra. A Tribute to the Work of Robert Gilmer'',
Brewer, J.W., Glaz, S., Heinzer, W.J., Olberding, B.M. (Eds.), Springer, 2006,  pp. 243-261. 

\bibitem{hmm} E. Houston, S. Malik, and J. Mott, {\em
Characterizations of $*$-multiplication domains}, Canad. Math.
Bull. 27(1984), 48-52.

%


\bibitem{hz} E. Houston and M. Zafrullah, {\em On $t$-invertibility}, II,
Comm. Algebra 17(1989), 1955-1969.

\bibitem{hz2} E. Houston and M. Zafrullah, {\em UMV-domains},
 in ``Arithmetical Properties of Commutative Rings and Monoids'',
 Lecture Notes   Pure  Appl. Math., Chapman and Hall,
  241(2005), 304-315.


%
  



%




%

  

   


  
\bibitem{mz}  J.L. Mott, and M. Zafrullah, {\em On Pr\"ufer $v$-multiplication domains},  Manuscripta Math. 35 (1981), 1-26.

  

%

  
\bibitem{o-m} A. Okabe and R. Matsuda, {\em Semistar operations on integral domains}, Math. J. Toyama Univ. 17(1994), 1-21.
  


\bibitem{n}
 M. Nagata, Local rings,  New York, Interscience, 1962.

%


\bibitem{pi} G. Picozza, Semistar operations and multiplicative ideal theory, Ph.D. Thesis, Universit\`a degli Studi ``Roma Tre'', 2004.


%

\bibitem{p} N. Popescu, Abelian categories with applications to rings and modules, Academic Press, New York, 1973.

\bibitem{st}
B. Stenstr\"om, Rings of quotients, Springer, New York, 1975.

%
%
%

%



%

  

\bibitem{WMc97} F.G. Wang and R.L. MacCasland, {\em On $w$-modules over strong Mori domains,} Comm. Algebra 25(1997), 1285--1306.

%
   
   

\end{thebibliography}

\end{document}